\documentclass{article}
\usepackage{amsthm}
\usepackage{graphicx} 
\usepackage{authblk}
\usepackage{longtable}
\usepackage[a4paper,margin=3cm]{geometry}
\usepackage{booktabs,multirow,tabularx} 
\usepackage{multicol} 
\usepackage{float}
\usepackage{amsmath, amssymb, mathtools}
\usepackage{supertabular,booktabs}
\usepackage{makecell}
\usepackage{caption}
\usepackage{setspace}
\usepackage{tikz}
\usepackage{qtree}
\usepackage{tikz-qtree}
\usepackage{subfigure}
\usepackage{algpseudocode}
\usepackage{algorithm}
\usepackage{multirow}

\newtheorem{theorem}{Theorem}[section]

\algnewcommand\algorithmicforeach{\textbf{for each}}
\algdef{S}[FOR]{ForEach}[1]{\algorithmicforeach\ #1\ \algorithmicdo}

\title{The Largest Condorcet Domain on 8 Alternatives}

\author[1]{Charles Leedham-Green}
\author[2]{Klas Markstr{\"o}m}
\author[3]{S\o ren Riis\thanks{Corresponding author}}
\affil[2]{University of Umeå}
\affil[1,3]{Queen Mary University of London}
\date{}
\begin{document}
\maketitle

\begin{abstract}
In this note, we report on a Condorcet domain of record-breaking size for $n$=8 alternatives. We show that there exists a Condorcet domain of size 224 and that this is the largest possible size for 8 alternatives. Our search also shows that this  domain is unique up to isomorphism and reversal. In this note we investigate properties of the new domain and relate them to various open problems and conjectures.

\end{abstract}

\section{Introduction}
\label{sec:intro}
Condorcet domains (CD), which are sets of linear orders giving rise to voting profiles with an acyclic pairwise majority relation, have  been studied by mathematicians, economists, and mathematical social scientists since the 1950s \cite{black, Arrow1951}. Condorcet domains find use in Arrovian aggregation and social choice theory \cite{puppe2019condorcet, Lackner17}. In social choice theory, a Condorcet winner is a candidate who would win over every other candidate in a pairwise comparison by securing the majority of votes \cite{monjardet2005social}. However, the existence of such a candidate is not always guaranteed, leading to the relevance of Condorcet Domains. 
A central question in this field has revolved around identifying large Condorcet domains, see    Fishburn, Gamlambos \& Reiner, Monjardet, Danilov \& Karzanov \& Koshevoy, Puppe \& Slinko, Karpov \& Slinko, Karpov 
\cite{fishburn1997acyclic,galambos2008acyclic, monjardet2009acyclic, danilov2012condorcet, puppe2022maximal,karpov2022constructing, karpov2022structured}.

 A significant category of Condorcet domains is rooted in Fishburn's alternating scheme, which alternates between two restriction rules on a subset of candidates and has been employed to construct numerous maximum size Condorcet domains. We refer to such domains based on the alternating scheme as Fishburn domains.
 
 Fishburn introduced a function $f(n)$ in \cite{fishburn1997acyclic}, defined to be the maximum size of a Condorcet domain on a set of $n$ alternatives, and posed the problem of determining the growth rate for $f(n)$. Fishburn also proved that for $n=16$ the Fishburn domain is not the largest CD.  This was followed by further research and bounds on $f(n)$ by Gamlambos \& Reiner, Danilov \& Karzanov, and Monjardet \cite{galambos2008acyclic,danilov2012condorcet,monjardet2009acyclic}. Karpov \& Slinko extended and refined this work in \cite{karpov2022symmetric} and Zhou \& Riis \cite{zhou2023new}.

Although extensive research has been conducted, all known maximum-sized Condorcet domains have been built using components based on either Fishburn's alternating schemes or his replacement scheme. For instance, Karpov\& Slinko \cite{karpov2022constructing} introduced a novel construction that enabled the creation of new Condorcet domains with unprecedented sizes. This allowed the authors to construct a Condorcet domain, superseding the size of Fishburn's domain for 13 alternatives. Recently, Zhou \& Riis \cite{zhou2023new} constructed Condorcet domains on 10 and 11 alternatives, superseding the size of the corresponding Fishburn domains. 

This paper shows that $n=8$ is the smallest number of  alternatives for which the Fishburn domain (size 222) is not the largest and that there is a Condorcet domain of size 224. Furthermore, relying on extensive computer calculation on the super-computer Abisko at Umeå, we also established 224 as an upper bound and that there, up to isomorphism, is only one such Condorcet domain. The need for a supercomputer, and a carefully devised algorithm, reflects the fact that a naive search would lead to search-tree with more than $6^{112}$ vertices. We also analyse some of the properties of this new domain.

\section{Preliminaries}
\label{sec:preliminaries}

There are many equivalent definitions of Condorcet domains.
In this paper, we adopt the definition proposed by Ward in~\cite{ward65}. According to this definition, a Condorcet domain of degree $n\geq 3$ is a set of orderings of $X_n=\{1,2,\ldots,n\}$ that satisfies certain local conditions. 

Specifically, a Condorcet domain of degree $n=3$ is defined as a set of orderings of $X_3$ that satisfies one of nine laws, denoted by $x$N$i$, where $x$ is an element of $X_3$, and $i$ is an integer between 1 and 3. The law $x$N$i$ requires that $x$ does not come in the $i$-th position in any order in the Condorcet domain. For example, $x$N$1$ means that $x$ may never come first, while $x$N$3$ means that $x$ may never come last.

A Condorcet domain of degree $n>3$ is a set $A$ of orderings of $X_n$ that satisfies the following property: the restriction of $A$ to every subset of $X_n$ of size 3 is a Condorcet domain. In other words, for every triple ${a,b,c}$ of elements of $X_n$, one of the nine laws $x$N$i$ must be satisfied, where $x\in{a,b,c}$. For example, $c$N$2$ would mean that $c$ may not come between $a$ and $b$ in any orderings in $A$.

A \emph{maximal Condorcet domain} of degree $n$ is a Condorcet domain of degree $n$ that is maximal under inclusion among the set of all Condorcet domains of degree $n$. A Maximum Condorcet domain is a Condorcet domain of the largest possible size for a given value of $n$.

To avoid repetition, we will use the acronyms CD and MCD,  to refer to Condorcet domain and Maximal Condorcet domain respectively.

For the case of degree 3, there are nine MCDs, each corresponding to one of the nine different laws $x$N$i$. It is easy to verify that these nine MCDs contain exactly four elements: two  transpositions and two even permutations (either the identity or a 3-cycle). Among the 9 MCDs of order 3, precisely six contain the identity order $1>2>3$ since the laws 1N1, 2N2, and 3N3 each rule out one CD of degree 3.

\subsection{Transformations and isomorphism of Condorcet domains}
First, recall that each linear order $A$ in a CD $B$ can also be viewed as a finite sequence of integers, obtained by ordering the elements of $X_n$ according to the linear order, and as the  permutation which permutes $X_n$ to this sequence.  We let $S_n$ denote the set of all permutations on $X_n$.

Let $g\in S_n$ and $i\in X_n$.  We define $ig$ as $g(i)$; and if $A$ is a sequence of elements of $X_n$ we define $Ag$ to be the sequence obtained by applying $g$ to the elements of $A$ in turn.  If $B$ is a CD, regarded as a set of sequences, we define $Bg$ to be the set of sequences obtained by applying $g$ to the sequences in $B$, and then $Bg$ is also a CD. Specifically, if $B$ satisfies the law $x$N$i$ on a triple $(a,b,c)$ for some $x\in(a,b,c)$, then $Bg$ satisfies the law $xg$N$i$ on the triple $(ag,bg,cg)$. We call CDs $B$ and $Bg$ \emph{isomorphic}. Therefore, two isomorphic CDs differ only by a relabelling of the elements of $X_n$.

The \emph{core} of a CD $B$ is the set of permutations $g\in B$ such that $Bg=B$. The core of a CD which contains the identity permutation $B$ is a group. We provide a more detailed discussion of the core in \cite{n7paper}.

\begin{table}[H]
    \renewcommand{\arraystretch}{1.3}
    \centering
    \begin{tabular}{ccccc}
    \toprule
   \textbf{Triple} & \textbf{Rule assigned} & \textbf{Condorcet domains} &  & \textbf{Core}\\

    \midrule
    \multirow{6}{*}{(i, j, k)}
    & 1N3 &
    \multirow{2}{*}{
    $\begin{rcases*}
        \begin{tabular}{cccc}
             ijk & jik & ikj & kij   \\
             ijk & jik & jki & kji \\
        \end{tabular}  
     \end{rcases*}$ 
    } & \multirow{2}{*}{isomorphic} & $\{ijk,ikj\}$  \\
    & 2N3 & && $\{(ijk),(kji)\}$ \\

    \cmidrule(l){2-5}

    & 3N1 &
    \multirow{2}{*}{
    $\begin{rcases*}
        \begin{tabular}{cccc}
             ijk & ikj & jik & jki  \\
             ijk & ikj & kij & kji  \\
        \end{tabular}  
     \end{rcases*}$ 
    } & \multirow{2}{*}{isomorphic}& $\{ijk,jik\}$ \\
    & 2N1 & && $\{(ijk),(kji)\}$\\

    \cmidrule(l){2-5}

    & 1N2 &
    \multirow{2}{*}{
    $\begin{rcases*}
        \begin{tabular}{cccc}
             ijk & ikj & jki & kji \\
             ijk & jik & kij & kji\\
        \end{tabular}  
     \end{rcases*}$ 
    } & \multirow{2}{*}{isomorphic}& $\{ijk,ikj\}$ \\
    & 3N2 & && $\{(ijk),(jik)\}$\\
    
    \bottomrule
    
\end{tabular}
    \caption{The Condorcet domains for 3 alternatives which contain the identity order. Each rule assigned to the triplet (i, j, k) with i$<$j$<$k is associated with a CD (which is given on the same line). The CDs displayed fall into 3 isomorphism classes, and each CD has a core of size 2. }
    \label{tab:cd_3}
\end{table}

It can be readily shown that for any Condorcet domain, the total number of 1N3 and 2N3 rules remains invariant under isomorphism. Likewise, this holds for the total number of 2N1 and 3N1 rules and the total number of 1N2 and 3N2 rules.

\section{Search methodology}
We developed an algorithm to generate all MCDs of a given degree $n$ and size at least equal to a user-specified  cutoff value (e.g. size $\geq 222$ for $n=8$).   We implemented this algorithm in C in a serial version which is sufficient for $n\leq 6$, and a parallelized version that we used for $n=7$ and $8$.  It is important to stress that this algorithm, unlike the one used by Zhou \& Riis \cite{zhou2023new}, aims to construct all MCDs above some user-specified size.

Our algorithm works by starting with the unrestricted domain of all linear orders on $n$ alternatives and then stepwise applying never laws $i$N$p$ to those triples which do not already satisfy some such law. The algorithm works with \emph{unitary} CDs, meaning CDs which contain the identity permutation. Since  every CD is isomorphic to some unitary CD this is without loss of generality. However, by using unitary CDs we reduce the set of possible never laws from 9 to 6, thereby speeding up our search. We will next sketch some of the details required in order to see that the algorithm is complete, though at first inefficient, and then how to also make it efficient.

We define the {\it Condorcet tree} of rank $n$,  which is a   homogeneous rooted tree of valency 6 and depth $ {n\choose 3}$,  as follows. The ${n \choose 3}$ triples of elements of $X_n$ are arranged in some order, so that the vertices of the tree at a given depth $t$ are associated with the corresponding triple $T_t$.  The six laws that a unitary CD may obey on a given triple are also ordered, and each child $w$ of a non-leaf $v$ of the tree
is associated with one such law $L_w$.  Every vertex $v$ is associated with a set  $c_v$ of linear orders on $X_n$.  If $v$ is the root then $c_v$ is the set of all orderings.
If $w$ is a child of $v$, where $v$ has depth $t$, then $c_w$ is obtained from $c_v$ by removing those orderings that do not satisfy the law $L_w$ when applied to $T_t$.

It is possible, in theory, to process the entire tree, depth first, constructing the sets $c_v$ for every vertex $v$.  Then the unitary  MCDs of degree $n$, as well as many non-maximal CDs, are found among the sets $c_v$ for the leaves $v$. In practice this is impracticable for $n>5$ as the tree is too big.  

For any leaf $v$ the set $c_v$ is a unitary CD, but these are not always maximal, and there will be very many duplicates. This arises from the fact that, as we move down the tree, the sets $c_v$ will often not only obey the laws that have been explicitly applied on triples but may also obey laws on triples which are implied by the applied laws. Using this observation allows a massive reduction in the number of vertices that need to be processed, giving us a tree with 0, 1 or 6 descendants from $v$ depending on whether $c_v$ cannot be maximal or must be a duplicate, has an implied law, or is unrestricted by earlier laws. This is determined as follows. 

When a vertex $v$ of height $t$ is processed  the law that was enforced on each triple $T_s$ for $s\le t$  to define $v$ - in other words the path from the root to $v$ - is recorded, and $c_v$ is constructed by taking $c_u$, where $u$ is the parent of $v$, and deleting all elements that do not satisfy the corresponding never  law $N_v$. For each $s\le t+1$ the set $L_s$ of laws that all the elements of $c_v$ obey when applied to the triple $T_s$ is determined. If, for some $s\le t$, the set $L_s$ contains a law that precedes the law $N_u$, where $u$ is the ancestor of $v$ of depth $s$, then the vertex $v$ is not processed any further, on the grounds of duplication, and its descendants are not visited.  Otherwise, for each $s\le t$, a law from $L_s$ is selected, and the set of sequences that obey all these laws is computed.  This set clearly contains $c_v$, and if, for some such selection of laws, this set strictly contains $c_v$ then again $c_v$ is not processed further.  In this case, any unitary CD arising from a leaf descendant of $v$ must either fail to be maximal, or will be a duplicate of a unitary MCD constructed from a descendant of another vertex of depth $t$.  If $v$ passes these tests, and $L_{t+1}$ is non-empty, the only descendant of $v$ that will be processed is the child $w$ defined by the least element of $L_{t+1}$, and then $c_w=c_v$.  Otherwise all children of $v$ are processed.

The validity of these restrictions of the full Condorcet tree  follows from a recursive argument which is given in full in \cite{n7paper}.

\section{Condorcet domains on 8 alternatives with size 224}
\label{sec:10}
Relying on extensive computer calculation on the super-computer Abisko at Umeå, we have  established that: 

\begin{theorem}
The maximum size of a CD on 8 alternatives is 224. Up to isomorphism, there is only one such CD. This CD has a core of size $4$. There are no MCDs of size 223.

The largest Condorcet domain containing the identity permutation and its reverse for $n=8$ alternatives is the Fishburn domain, which has a size of 222. 
\end{theorem}
We aim to extend this with more precise counts and analysis of other large Condorcet domains on 8 alternatives in an upcoming paper. 

Now let us investigate the properties of the MCD of size 224.

\begin{enumerate}
    \item The Fishburn domain has size 222 and hence is not the maximum CD for $n=8$  alternatives
    
    \item There are 56 isomorphic Condorcet domains of size 224 which contain the identity order. Among these there is one special MCD we will refer to as D224, where each never-rule - except for the two triplets (123) and (678) - is 1N3 or 3N1. We display the rules for D224 in Table \ref{table:triplets_rules_8} and its linear orders in  Table \ref{table:triplets_rules_8B}
    
    \item  The domain does not have maximal width, i.e. it does not contain a pair of reversed orders.
    
    \item The domain is self-dual. That is, the domain is  isomorphic to the domain obtained by reversing each of its linear orders. 

    \item The restriction of the domain to each triple of alternatives has size 4. This means that this domain is copious in the terminology of \cite{SLINKO2019166} and is equivalent to the fact that the domain satisfies exactly one never-rule on each triple. 

    \item  The domain is a peak-pit domain in the sense of \cite{danilov2012condorcet}, i.e. every triple satisfies a condition of either the form $x$N1 or $x$N3, for some $x$ in the triple.  

    \item The authors of \cite{karpov2022constructing} asked for examples of maximum CDs which are not peak-pit domains of maximal width. Our domain is the first known such example and shows that $n=8$ is the smallest $n$ for which this occurs.

    \item The domain is connected (see \cite{monjardet2009acyclic} for the lengthy definition of this well used property.)  This is in line with the conjecture from \cite{puppe2022maximal} that all maximal peak-pit CDs  are connected. 

    \item  The domain has a core of size 4, which is given in captions of Tables \ref{table:triplets_rules_8} and \ref{table:triplets_rules_8B}.

\end{enumerate}

\begin{table}[H]
\centering
\begin{tabular}{lc}
\toprule
\textbf{Triplets} & \textbf{Rules}\\
\midrule
(1, 2, 3) & 2N3 \\
(1, 2, 4) & 1N3 \\
(1, 2, 5) & 3N1 \\
(1, 2, 6) & 3N1 \\
(1, 2, 7) & 3N1 \\
(1, 2, 8) & 3N1 \\
(1, 3, 4) & 1N3 \\
(1, 3, 5) & 3N1 \\
(1, 3, 6) & 3N1 \\
(1, 3, 7) & 3N1 \\
(1, 3, 8) & 3N1 \\
(1, 4, 5) & 1N3 \\
(1, 4, 6) & 1N3 \\
(1, 4, 7) & 1N3 \\

\bottomrule
\end{tabular}
\hspace{1em}
\begin{tabular}{lc}
\toprule
\textbf{Triplets} & \textbf{Rules}\\
\midrule
(1, 4, 8) & 1N3 \\
(1, 5, 6) & 1N3 \\
(1, 5, 7) & 1N3 \\
(1, 5, 8) & 3N1 \\
(1, 6, 7) & 3N1 \\
(1, 6, 8) & 1N3 \\
(1, 7, 8) & 1N3 \\
(2, 3, 4) & 1N3 \\
(2, 3, 5) & 1N3 \\
(2, 3, 6) & 1N3 \\
(2, 3, 7) & 1N3 \\
(2, 3, 8) & 1N3 \\
(2, 4, 5) & 3N1 \\
(2, 4, 6) & 3N1 \\

\bottomrule
\end{tabular}
\hspace{1em}
\begin{tabular}{lc}
\toprule
\textbf{Triplets} & \textbf{Rules}\\
\midrule
(2, 4, 7) & 3N1 \\
(2, 4, 8) & 3N1 \\
(2, 5, 6) & 1N3 \\
(2, 5, 7) & 1N3 \\
(2, 5, 8) & 3N1 \\
(2, 6, 7) & 3N1 \\
(2, 6, 8) & 1N3 \\
(2, 7, 8) & 1N3 \\
(3, 4, 5) & 3N1 \\
(3, 4, 6) & 3N1 \\
(3, 4, 7) & 3N1 \\
(3, 4, 8) & 3N1 \\
(3, 5, 6) & 1N3 \\
(3, 5, 7) & 1N3 \\
\bottomrule
\end{tabular}
\hspace{1em}
\begin{tabular}{lc}
\toprule
\textbf{Triplets} & \textbf{Rules}\\
\midrule

(3, 5, 8) & 3N1 \\
(3, 6, 7) & 3N1 \\
(3, 6, 8) & 1N3 \\
(3, 7, 8) & 1N3 \\
(4, 5, 6) & 1N3 \\
(4, 5, 7) & 1N3 \\
(4, 5, 8) & 3N1 \\
(4, 6, 7) & 3N1 \\
(4, 6, 8) & 1N3 \\
(4, 7, 8) & 1N3 \\
(5, 6, 7) & 3N1 \\
(5, 6, 8) & 3N1 \\
(5, 7, 8) & 3N1 \\
(6, 7, 8) & 2N1 \\
\bottomrule
\end{tabular}

\caption{Table of triplets and rules that produces the Condorcet domain D224 of size 224 for 8 alternatives. This specific CD is invariant under
the action by the permutations group $G =\{{\rm id}, (12)(34), (56)(78), (12)(34)(56)(78)\}$}
\label{table:triplets_rules_8}
\end{table}

\begin{table}[H]
  \hrule
  \vspace{1ex}
  \begin{center}
    \textbf{Condorcet domain with 224 Permutations for 8 Alternatives}    
  \end{center}
  \vspace{-2ex}
  \hrule 
  \vspace{2ex}
\begin{spacing}{1}
\noindent 
\underline{12345678}
12345687 
12345867 
12345876 
12346578 
\underline{12346587} 
12346758 
12346785 
12354678 
12354687 
12354867 
12354876  
12358467  
12358476 
12364578 
12364587 
12364758 
12364785 
12367458 
12367485 
12435678 
12435687 
12435867 
12435876 
12436578 
12436587 
12436758 
12436785 
12453678 
12453687 
12453867 
12453876 
12458367 
12458376 
12463578 
12463587 
12463758 
12463785 
12467358 
12467385 
14235678
14235687
14235867
14235876
14236578
14236587
14236758
14236785
14253678
14253687
14253867
14253876
14258367
14258376
14263578
14263587
14263758
14263785
14267358
14267385
14523678
14523687
14523867
14523876
14528367
14528376
14582367
14582376
14623578
14623587
14623758
14623785
14627358
14627385
14672358
14672385
21345678
21345687
21345867
21345876
21346578
21346587
21346758
21346785
21354678
21354687
21354867
21354876
21358467
21358476
21364578
21364587
21364758
21364785
21367458
21367485
\underline{21435678}
21435687
21435867
21435876
21436578
\underline{21436587}
21436758
21436785
21453678
21453687
21453867
21453876
21458367
21458376
21463578
21463587
21463758
21463785
21467358
21467385
23145678
23145687
23145867
23145876
23146578
23146587
23146758
23146785
23154678
23154687
23154867
23154876
23158467
23158476
23164578
23164587
23164758
23164785
23167458
23167485
23514678
23514687
23514867
23514876
23518467
23518476
23581467
23581476
23614578
23614587
23614758
23614785
23617458
23617485
23671458
23671485
32145678
32145687
32145867
32145876
32146578
32146587
32146758
32146785
32154678
32154687
32154867
32154876
32158467
32158476
32164578
32164587
32164758
32164785
32167458
32167485
32514678
32514687
32514867
32514876
32518467
32518476
32581467
32581476
32614578
32614587
32614758
32614785
32617458
32617485
32671458
32671485
41235678
41235687
41235867
41235876
41236578
41236587
41236758
41236785
41253678
41253687
41253867
41253876
41258367
41258376
41263578
41263587
41263758
41263785
41267358
41267385
41523678
41523687
41523867
41523876
41528367
41528376
41582367
41582376
41623578
41623587
41623758
41623785
41627358
41627385
41672358
41672385
\end{spacing}
\caption{Permutation in Condorcet domain corresponding to the rules in table \ref{table:triplets_rules_8}} The CD's core consists of the underlined permutations 12345678, 12346587,21435678 and 21436587.
\label{table:triplets_rules_8B}
\end{table}

\section{Conclusion} In conclusion, our work has demonstrated a record-breaking maximum Condorcet domain 
for $n=8$ alternatives,  which is essentially unique (up to isomorphism and reversal). We have also investigated how our domain relates to various well-studied properties of MCDs.   Our findings contribute to understanding the structure of
Condorcet domains and have potential applications in voting theory and social choice. 

Overall, our work highlights the importance of understanding the properties and structures of CDs in order to construct larger examples  and might pave the way for future research in this area.

We also observe that some record-breaking CDs for $n=8$ alternatives exhibit almost all rules of the form 1N3 and 3N1. These rules can be interpreted as a form of seeded voting. In such a system, for each set of three alternatives, a seeding is implemented to restrict the lowest-seeded alternative from being the highest-ranked preference or the highest-seeded alternative from being the lowest-ranked preference. A better understanding  of the global effects of this type of local seeding could serve as a foundation for future research, potentially offering insights into algorithmic fairness and impartiality in computer-supported decision-making.

\section*{Acknowledgements}
This research was conducted using the resources of High Performance Computing Center North (HPC2N). We would like to thank the anonymous reviewers for their constructive criticism.

\section*{Statements and Declarations}
The authors are listed alphabetically and declare no conflict of interest. 

\bibliographystyle{plain}
\bibliography{references}

\end{document}